# Charles Babbage, Ada Lovelace, and the Bernoulli Numbers

Thomas J. Misa

University of Minnesota

ABSTRACT: This chapter assembles the pertinent sources to suggest several corrections to an unduly negative scholarly view of Ada Lovelace.  I suggest that the Lovelace-Babbage question is not a zero-sum game, where any credit added to Lovelace somehow detracts from Babbage.  Ample evidence indicates Babbage and Lovelace each had important contributions to the 1843 Sketch and the accompanying Notes.  Further, prior claims about her lack of mathematical background seem doubtful after consulting Lovelace's detailed correspondence with two highly accomplished figures in 19th century mathematics, Charles Babbage and Augustus De Morgan.  Babbage and Lovelace's treatment of the Bernoulli numbers in note "G" spotlights the mathematical sophistication their collaboration.  Finally, acknowledging significant definitional problems in calling Lovelace the world's "first computer programmer," I conclude that Lovelace created an step-by-step elemental sequence of instructions—that is, an algorithm—for computing the series of Bernoulli numbers that was intended for Babbage's Analytical Engine.

Few figures in the long history of computing generate more passion and sometimes more enmity than Charles Babbage and Ada Lovelace.  History has treated Babbage as a brilliant but temperamental pioneer in a half dozen scientific fields, an "irascible genius" in one biographer's persisting image.  Histories of mathematics typically praise his efforts to bring the modern notation of continental calculus to Cambridge University where the long shadow of Isaac Newton had reigned for more than a century. Babbage was made a Fellow of the Royal Society in 1816, just two years after leaving Cambridge.  In 1820, he founded the Astronomical Society to standardize observational data and improve positional calculations.  Babbage immersed himself in numerous projects and publications during the next two decades.  While histories of insurance credit a early book on actuarial calculations (1826), histories of science debate his polemical *On the Decline of Science* (1830), histories of industry spotlight *On the Economy of Machinery and Manufactures* (1832), and histories of religion and science note his *Ninth Bridgewater Treatise* (1837).  In the midst of this publishing storm, he served as Lucasian Professor of Mathematics at Cambridge (1828-39), Newton's old post; inherited a sizable fortune from his father; married, raised a family, and lost his wife; and embarked upon designing two mechanical computing machines, the simple but elegant Difference Engine and the complex and enigmatic Analytical Engine.  During these years he hosted a fashionable salon gathering in London, twice ran for Parliament, traveled widely, and corresponded energetically.  He made some powerful friends and a few powerful enemies.

In the early 1840s Babbage collaborated closely with Ada Lovelace, and this essay examines their work as an intellectually intimate "creative couple."[1]  Especially in the popular mind, Ada Lovelace has recently been on a roll.  She is lionized as the founder of scientific computing and hailed as the world's first computer programmer.  "Readers will recognize Steve Jobs, Charles Babbage, Bill Gates, and Ada Lovelace as appropriate inclusions in [the young-adult book] Computer Technology Innovators," states a 2013 review in *School Library Journal*.[2]  And in his recent best-selling *The Innovators*, Walter Isaacson employs Ada Lovelace as bookends: in chapter 1, "Ada, Countess of Lovelace," it is she (above Babbage) who is an engaging founding figure; and in his concluding chapter, "Ada Forever," he places the promise of computer innovation in the hands of her "spiritual heirs."  As chapters in this volume amply attest, Ada's legacy is wide and deep.  There are no other 19th century women who have a programming language named for them (chapters 3-5) and figure prominently in a contemporary science-fiction literary genre (chapters 8-10) and serve as inspiration for contemporary computing reform (chapters 11-13).

Scholarship on these two figures is a something of a puzzlement.  A scientific figure like Babbage should have inspired a full-length biography somewhere along the line, but despite numerous essays and several books, we still lack a complete life-and-times biography.[3]  One recent effort by David Alan Grier to explore such a biography confronted a daunting mass of archival materials, some in private hands and difficult to access.  Babbage's published papers alone run to eleven volumes.[4]  Several popular treatments of Ada Lovelace have appeared, in addition to Isaacson's,[5] but the existing scholarly consensus on her is a dour one, often highly critical. Allan Bromley, whose research in the Babbage materials inspired Doron Swade and led to the Science Museum's project to reconstruct the Difference Engine No. 2, saw the world from a Babbage-centered perspective. "The [Lovelace] translation has extensive notes, *written under Babbage's supervision*, that give an excellent account of *his* understanding of the mechanization of computational processes and the mathematical powers of the machine," he wrote in 1982 (emphasis added). Bromley's dismissal of Lovelace hardened in subsequent years.[6] Dorothy Stein, in *Ada: A Life and a Legacy* (1985), sternly cautioned that Lovelace was "a figure whose achievement turns out not to deserve the recognition accorded it."[7] Stein's and Bromley's gloomy verdict persists in the scholarly survey *Computer: History of the Information Machine* (1996), now in its third edition (2014), where the key critical passage on Lovelace remains: "the extent of Lovelace's intellectual contribution to the *Sketch* has been much exaggerated . . . . Later scholarship has shown that most of the technical content and all of the programs in the *Sketch* were Babbage's work."[8] The most strident negative verdict derives from Bruce Collier's 1970 Harvard thesis, recently publicized in the *Economist* magazine on Ada Lovelace day:

> Ada was as mad as a hatter, and contributed little more to the "Notes" than trouble . . . . I will retain an open mind on whether Ada was crazy because of her substance abuse . . . or despite it. I hope nobody feels compelled to write another book on the subject. But, then, I guess *someone* has to be the most overrated figure in the history of computing.[9]

---

[1] Helena M. Pycior, Nancy G. Slack, and Pnina G. Abir-Am., eds, *Creative Couples in the Sciences* (New Brunswick, NJ: Rutgers University Press, 1996).  "There isn't a hint of romance in any of their correspondence with one another," according to Sydney Padua's *The Thrilling Adventures of Lovelace and Babbage* (New York: Pantheon, 2015), quote 38 note 16.  Babbage, a widower, was the age of Ada's mother.

[2] Vicki Reutter, "Computer Technology Innovators," *School Library Journal* (Oct. 2013): 65.

[3] Anthony Hyman's *Charles Babbage: Pioneer of the Computer* (Princeton: Princeton University Press, 1982) treats Babbage's life in 250 pages.  An early study was Maboth Moseley's *Irascible Genius: A Life of Charles Babbage, Inventor* (London: Hutchinson, 1964), which is attacked by Dorothy Stein, in *Ada: A Life and a Legacy* (Cambridge: MIT Press, 1985), p. x, as "almost perversely inaccurate, distorted, and fabricated."  A specialized and valuable study is J. M. Dubbey, *The Mathematical Work of Charles Babbage* (Cambridge: Cambridge University Press, 1978). *[See also Christopher Hollings, Ursula Martin, and Adrian Rice's* Ada Lovelace: The Making of a Computer Scientist *(Bodleian Libraries, 2018)*

[4] David Alan Grier, "The Inconsistent Youth of Charles Babbage," *IEEE Annals of the History of Computing* 32 no. 4 (2010): 18-31; Martin Campbell Kelly, ed., *The Works of Charles Babbage* (London: Pickering / New York: New York University Press, 1989; 11 vols.).

---

[5] James Essinger, *A Female Genius: How Ada Lovelace Started the Computer Age* (London: Gibson Square, 2013), which appeared in the United States as *Ada's Algorithm* (Brooklyn: Melville House, 2014). Essinger had earlier written *Jacquard's Web: How a Hand-loom Led to the Birth of the Information Age*(Oxford: Oxford University Press, 2004). Contrast Martin Davis and Virginia Davis, "Mistaken Ancestry: The Jacquard and the Computer," *Textile* 3 no. 1 (2005): 76-87. Earlier positive treatments include Betty Alexandra Toole, ed., *Ada, The Enchantress of Numbers: A Selection from the Letters of Lord Byron's Daughter* (Mill Valley CA: Strawberry Press, 1992) and Doris Langley Moore, *Ada, Countess of Lovelace: Byron's Legitimate Daughter* (New York: Harper & Row, 1977).

[6] Doron D. Swade, "Redeeming Charles Babbage's Mechanical Computer," *Scientific American* (February 1993): 86-91; Allan Bromley, "Charles Babbage's Analytical Engine, 1838," *Annals of the History of Computing* 4 no. 3 (1982): 196-217, quote p. 197; Allan Bromley, "Difference and Analytical Engines," in William Aspray, ed., *Computing before Computers* (Ames: Iowa State University Press, 1990), 59-98, on p. 89.

[7] Dorothy Stein, *Ada: A Life and a Legacy* (Cambridge: MIT Press, 1985), quote xii. Walter Isaacson's *The Innovators: How a Group of Hackers, Geniuses, and Geeks Created the Digital Revolution* (New York: Simon & Schuster, 2014), on p. 493 note 1, praises Stein's book as "the most scholarly and balanced." A later negative view is Jay Belanger and Dorothy Stein, "Shadowy Vision: Spanners in the Mechanization of Mathematics," *Historia Mathematica* 32 (2005): 76-93. Stein strongly criticized Dorothy L. Moore, *Ada, Countess of Lovelace: Bryon's Legitimate Daughter* (London: Murray, 1977).

[8] Martin Campbell-Kelly, William Aspray, Nathan Ensmenger, and Jeffrey R. Yost, *Computer: A History of the Information Machine* (Boulder CO: Westview Press, 2014), quote p. 44. Note that Bromley qualified his claim: that all but one program (Bernoulli numbers) was Babbage's. Compare Campbell Kelly and Aspray's first edition of *Computer* (p. 57).

In this essay, I assemble the pertinent sources, including correspondence about the Notes to the Menabrea Sketch (printed in this volume: LINK), and suggest several modest corrections to the unduly negative scholarly view.[10] First, in contrast to much of the existing literature, the Lovelace-Babbage question is not a zero-sum game, where some portion of credit added to Lovelace somehow detracts from Babbage, or vice versa. There is ample evidence that Babbage and Lovelace each had important contributions to the Sketch and the Notes, and attention to their intellectual collaboration is revealing. Second, claims about her lack of mathematical background seem doubtful after consulting Lovelace's detailed correspondence with Babbage and Augustus De Morgan, two highly accomplished figures in 19th century mathematics. The treatment of the Bernoulli numbers in note "G" spotlights the intellectually intimate collaboration between Babbage and Lovelace and its mathematical sophistication. Finally, while there may be significant definitional problems in calling Lovelace the world's "first computer programmer," the evidence is reasonably clear that Lovelace created an step-by-step elemental sequence of instructions—that is, an algorithm—for computing the series of Bernoulli numbers that was intended for Babbage's Analytical Engine. The underlying mathematics might well have been Babbage's, for he was a distinguished mathematical and scientific figure. Lovelace transformed an equation for the Bernoulli numbers into a precise series of elemental additions, multiplications, and substitutions.

The algorithm specified a sequence of calculations, requiring a real-life computer capable of running a program with a looping structure and conditional testing. Contemporary computer experts have noted in the large table (representing the Bernoulli number algorithm) in the Sketch's Note "G" there was one misplaced minus sign which, when corrected, led to the result that Ada's algorithm correctly computed the series of Bernoulli numbers. In "The Babbage Machine and the Origins of Programming," the authors reproduce Lovelace's table for the Bernoulli numbers and translate the algorithm into a 65-line FORTRAN program that computes them.  The program has eight "if . . . [then] go-to" statements and a simple structure, straight from Lovelace's table-algorithm, that builds up algebraic statements one mathematical operation at a time: for example, computing the expression (2n - 1) / (2n + 1) requires 4 program steps.[11]  In the original 1843 publication of Note "G," there are clearly two nested loops, embedded in a larger looping structure (see below).  There is direct documentary evidence that Ada Lovelace created this table (writing it out in pencil).  She and Babbage corresponded intensively in the weeks and days prior to its publication in *Taylor's Scientific Memoirs*.  This series, published in London between 1837 and 1852 by Richard Taylor in cooperation with the British Association for the Advancement of Science (BAAS), printed English translations of prominent European scientific papers.  "Here was to be found . . . such European leaders" as Bessel, Bunsen, Gauss, Ohm and many others.[12]  The Analytical Engine was Babbage's creation while the Sketch and Notes are best understood as the product of an intense intellectual collaboration between Babbage and Lovelace.

---

[9] "Ada Lovelace Day: Right Idea. Wrong Woman?" *Economist* (24 March 2010) at www.economist.com/node/21005551/print (accessed 2015). I use "derived from" intentionally, since in the printed copy of Collier's dissertation I examined (CBI QA75.C634x 1970a) I did not find this quotation; see also robroy.dyndns.info/collier (Jan. 2015). Doron Swade quotes Collier in *The Cogwheel Brain: Charles Babbage and the Quest to Build the First Computer* (London: Little, Brown, 2000), p. 168.

[10] Correspondence in the Toole, Stein, Moore, and Swade volumes as well as the document-centered accounts in Velma R. Huskey and Harry D. Huskey, "Lady Lovelace and Charles Babbage," *Annals of the History of Computing* 2 no. 4 (1980): 299-329; and John Fuegi and Jo Francis, "Lovelace & Babbage and the Creation of the 1843 'Notes'," *IEEE Annals of the History of Computing* 25 no. 4 (2003): 16-26.

**Babbage and Lovelace**

It is with good reason that computing history scholars have praised the research of Allan G. Bromley (1947-2002).  Bromley, a computer scientist at the University of Sydney, made careful studies of the Babbage letters and notebooks that led to the Science Museum's reconstruction of the Babbage Difference Engine No. 2.  In 2000 Tim Bergin, editor-in-chief of *IEEE Annals of the History of Computing*, introduced a special issue of the journal dedicated to Bromley's scholarship by noting his "fundamental contributions," former *Annals* editor-in-chief Michael Williams called his work "groundbreaking," and none other than computer pioneer Maurice Wilkes stated "it is to Bromley that we owe nearly all our present knowledge of Babbage's work on computing machinery at the detailed mechanical level." Wilkes noted Bromley's determined insistence that "Babbage's work on the Analytical Engine was completely original."[13]

---

[11] Campbell Kelly, *Works of Charles Babbage*, volume 3: 159 note a.  Garry J. Tee, in reviewing a 1979 Russian publication by A. K. Petrenko and O. L. Petrenko, wrote this: "The most advanced illustration given in Lovelace's paper is an elaborate program for computing the sequence of Bernoulli numbers, which was written by Babbage but corrected by her. The present authors have transcribed that program into FORTRAN, detecting thereby a few misprints but only one significant error, with one variable having the wrong sign. Their transcribed version is easier for a modern reader to understand than the original program of 1843, and it does compute correctly the sequence of Bernoulli numbers." See <www.ams.org/mathscinet-getitem?mr=83b:01045> (accessed 2015) reviewing Petrenko and Petrenko's "The Babbage Machine and the Origins of Programming," [in Russian] *Istoriko-matematicheskie issledovaniià* 24 (1979): 340-360, 389.  I examined the FORTRAN program in the Russian original.

[12] W. H. Brock and A. J. Meadows, *The Lamp of Learning: Taylor & Francis and Two Centuries Of Publishing* (London: Taylor & Francis, 1998; 2nd edition), esp. chapter 4 "Taylor and the Commercial Science Journal," quote p. 105.

By contrast, Bromley's view on Lovelace was strongly critical. "All but one of the programs cited in her notes had been prepared by Babbage from three to seven years earlier. The exception [on Bernoulli numbers] was prepared by Babbage for her, although she did detect a 'bug' in it," he wrote. "Not only is there no evidence that Ada Lovelace ever prepared a program for the Analytical Engine but her correspondence with Babbage shows that she did not have the knowledge to do so."[14] It is Bromley's viewpoint, slightly modified, that finds its way into the recent edition of the highly regarded *Computer: A History of the Information Machine* (2014) with its undue assertion that "all of the programs in the *Sketch* were Babbage's work." Reviewing the evidence about Lovelace's mathematical knowledge and the writing of the notes to her translation of Menabrea's sketch might modify these overly negative assertions. In later correspondence with Wilkes, Bromley did allow, in comparison with fellow Babbage scholar Doron Swade, "I have been known to express my views more intemperately."[15]

One must acknowledge that Ada Lovelace in her energetic, imaginative, self-absorbed, and at times grandiloquent correspondence gives her later-day critics much to aim at. Her self-regarding statements about her own mathematical abilities can be off-putting. And her pointed remarks sometimes aimed at Charles Babbage might rub the wrong way anyone who thinks distinguished scientists should be treated with dignity and respect.[16] Science at the time was expanding from its strictly male enclaves at Oxford and Cambridge, with the creation of the BAAS (f. 1831) and other learned societies, but in the new scientific institutions Ada Lovelace and Mary Somerville were treated at best as "second class members."[17] One need not make excuses for her imaginative flights of fancy, but it does need to be borne in mind that Lovelace was the daughter of an aristocratic Baron (the poet Lord Byron) and married to a highly ranked Earl. Her mother had schemes drawing on the family's network that extended into the royal family itself. Babbage, although he was handsomely wealthy after his banker-father died in 1827 and a fortune of £100,000 passed down to him, was all the same a commoner. Ada sought for years to land Babbage a knighthood.[18]

---

[13] Tim Bergin, "About this Issue" (quote p. 2); Michael R. Williams, "Allan Bromley," (quote p. 3); Maurice V. Wilkes, "Introduction to 'Babbage's Analytical Plans 28 and 28a—The Programmer's Interface'," (quote p. 4) in *IEEE Annals of the History of Computing* 22 no. 4 (2000).

[14] Allan Bromley, "Difference and Analytical Engines," in William Aspray, ed., *Computing before Computers* (Ames: Iowa State University Press, 1990), 59-98, quote p. 89.

[15] See "Appreciating Charles Babbage: Emails between Allan Bromley and Maurice Wilkes," *IEEE Annals of the History of Computing* 26 no. 4 (2004): 62-70, on 62 (intemperately).

[16] In the midst of writing the translation's notes, when letters and draft manuscripts were passed daily between them, Ada writes Babbage, somewhat curtly: "Now *pray* attend strictly to my requests; or you will cause me very serious annoyance," quoted in Huskey and Huskey, "Lady Lovelace and Charles Babbage," p. 313. There is a curious mix of defiance and deference in Ada's correspondence with Babbage, Somerville, De Morgan and other prominent figures.

[17] Ruth Watts, *Gender, Power and the Unitarians in England, 1760-1860* (New York: Longman, 1998), quote p. 153.

The first line of evidence suggesting an intellectual partnership between Charles Babbage and Ada Lovelace comes from witnesses to their first meeting. She first met Babbage in 1833, a year after being formally presented to the court, through an introduction by Mary Somerville, the mathematician, scientific popularizer, and notable English translator of Pierre-Simon Laplace's *Mécanique Céleste*; the two women kept up a scientific correspondence for many years.[19] Two years after she met Babbage, Ada married William King, already a Baron, who was within three years created the first Earl of Lovelace, and so it is a convenience to refer to her as Ada Lovelace rather than the more precise but ponderous Augusta Ada King, Right Honourable the Countess of Lovelace.

There is eyewitness evidence that, when she saw it, Ada grasped the principles and significance of Babbage's prototype Difference Engine. Babbage had started work on it in 1822, and it was in an advanced state of development in 1833; fatefully, the following year Babbage set the Difference Engine aside and focused instead on the conceptually elaborate Analytical Engine, which remained a "brilliant obsession" nearly to the end of his life.[20] In June 1833 Ada's mother, Lady Byron, described a visit along with her daughter and a friend to inspect Babbage's machine in some detail and with great wonder while admitting herself only "faint glimpses of the principles by which it worked." The Difference Engine was at the time able to compute polynomial expressions, extract roots to quadratic equations, and count to 10,000. Ada saw its significance. "I well remember accompanying her to see Mr. Babbage's wonderful analytical engine," wrote Sophia De Morgan, the wife of mathematician Augustus De Morgan and, like Mary Somerville, a long-term correspondent with Ada herself. "While other visitors gazed at the working of this beautiful instrument with the sort of expression . . . that some savages are said to have shown on first seeing a looking-glass or hearing a gun . . . Miss Byron, young as she was, understood its working and saw the great beauty of the invention. She had read the Differential Calculus to some extent, and after her marriage she pursued the study and translated a small work of the Italian mathematician Menabrea, in which the mathematical principles of its construction [were] explained."[21]

---

[18] Doron Swade, s.v. Babbage, Charles (1791-1871), *Oxford Dictionary of National Biography* (Oxford University Press, 2004); online edition, May 2009 at www.oxforddnb.com/view/article/962. Babbage died in 1871 with a fortune "under" £40,000.

[19] Elizabeth Chambers Patterson, *Mary Somerville and the Cultivation of Science, 1815-1840* (Boston: Nijhoff, 1983); and Kathryn A. Neeley, *Mary Somerville: Science, Illumination, and the Female Mind* (Cambridge: Cambridge University Press, 2001).

[20] See "Appreciating Charles Babbage: Emails between Allan Bromley and Maurice Wilkes," *IEEE Annals of the History of Computing* 26 no. 4 (2004): 62-70, quote p. 63 (obsession).

Babbage invited Ada with a friend or chaperone to attend his series of "Saturday evenings" where his London home became a fashionable salon, filled with celebrities from the political, cultural, and scientific world. Charles Dickens and Charles Darwin, among hundreds of others, were happy to mix with the well-cultured crowd. "One of three qualifications were necessary for those who sought to be invited—intellect, beauty, or rank—without one of these, you might be rich as Croesus—and yet be told, you cannot enter here," recalled one society figure. "His calculating machine was an endless subject of monologue."[22] Some were deeply impressed by its ability to generate a list of prime numbers. Since it could solve any second-degree polynomial equation, Babbage set it to compute a series of 40 prime numbers by evaluating the expression $x^2 + x + 41$ for the first 40 integers. A few months after the invitation, she wrote to Mary Somerville asking her to convey to Babbage's son "how exceedingly obliged I am . . . for his unexpected kindness in sending me the plates & account of the Machine, which is exactly what I was in want of; & is a very great help to me."[23] Ada was at the time 19 years old.

A significant line of evidence bearing on Bromley's claim that "she did not have the knowledge" to prepare a program for the Analytical Engine is the substantial depth of Ada's mathematical studies beginning in the 1830s, which predated her contact with the Menabrea manuscript and the translation project of the early 1840s. "Ada was much attached to me, and often came to stay with me. It was by my advice that she studied mathematics," recalled Mary Somerville. "She always wrote to me for an explanation when she met with any difficulty. Among my papers I lately found many of her notes, asking mathematical questions."[24]

During these years Ada had three tutors in mathematics, in addition to intellectual interchange with Somerville, Babbage, and her scientifically minded husband, who was made a Fellow of the Royal Society in 1841; two of them were distinguished figures. Her first mathematics tutor was the elderly William Frend, the notable social reformer who had authored

---

[21] Lady Byron quoted in Moore, *Ada*, p. 44; Sophia De Morgan, *Memoir of Augustus De Morgan* (London: Longmans, Green, 1882), quote p. 89 (accompanying her).

[22] "John Kenyon and His Friends," *Temple Bar: A London Magazine for Town and Country Readers* 88 (1890): quote p. 490 (three qualifications); AAL to Mary Somerville 19 March 1834 in Toole, *Ada*, p. 57 (Babbage's invitation).

[23] AAL to Mary Somerville 8 November 1834 quoted Huskey and Huskey, "Lady Lovelace and Charles Babbage," quote p. 303 (account of the Machine).

[24] Martha Somerville, ed., *Personal Recollections from early life to old age of Mary Somerville* (Boston: Roberts Brothers, 1874), quote p. 154. *[See also, two subsequent publications: Hollings et al "The Lovelace-De Morgan Mathematical Correspondence: A Critical Re-Appraisal,"* Historia Mathematica *2017 at https://doi.org/10.1016/j.hm. 2017.04.001 and "The early mathematical education of Ada Lovelace"* BSHM Bulletin *(2017) at https://doi.org/10.1080/17498430.2017.1325297*

*Principles of Algebra* (1796) along with many other tracts; among his students had been Ada's own mother.[25] It seems likely that Frend introduced Ada to Mary Somerville, connecting her with Babbage. Ada and Babbage corresponded on mathematical topics following their 1833 meeting, again years prior to the translation project. Perhaps her most important mathematical tutor was her friend Sophia's husband and William Frend's son-in-law, Augustus De Morgan. He gave Ada Lovelace, as Sophia later wrote, "much help in her mathematical studies, which were carried farther than her mother's had been."[26] Even the severely critical study by Dorothy Stein acknowledges the De Morgan-Lovelace letters as "a correspondence course in calculus."[27]

Augustus De Morgan was like Babbage a graduate of Cambridge, a disbeliever in the traditional Church of England, and a distinguished mathematician. He was named founding professor of mathematics at London University (now University College London), shortly after its founding in 1826, at the age of 22. It was a secular university, unlike Cambridge and Oxford, and admitted women as regular students. De Morgan published books on trigonometry, arithmetic, algebra, probability, and logic. In the early 1840s while exchanging regular letters with De Morgan on the topic, Ada reported that she was "drowning in Calculus." During these years De Morgan was working on a book project for the London-based Society for the Diffusion of Useful Knowledge (1826-48), published in 1842 as an 800-page textbook on *Differential and Integral Calculus*.[28] Her letters to De Morgan are filled with specific questions about differential calculus, limits, Leibnitz's notation, three-dimensional geometry, notation of functions, and standards of reasoning and proof. On 21 November 1841, in the context of her mathematical exercises, she asked him about the "numbers of Bernoulli."[29] Ada's awareness of the Bernoulli numbers thus *predated* her work on the Menabrea translation and the writing of Note G describing an algorithm for their computation. Whereas Ada's first tutor, the elder William

---

[25] Judith S. Lewis, "Princess of Parallelograms and Her Daughter: Math and Gender in the Nineteenth Century English Aristocracy," *Women's Studies International Forum* 18 no. 4 (1995): 387-394.

[26] Sophia De Morgan, *Memoir of Augustus De Morgan* (London: Longmans, Green, 1882), quote p. 89 (much help)

[27] Stein, *Ada*, quote pp. xii (correspondence course)

[28] Toole, Ada, quote p. 169 (drowning in Calculus). See Augustus De Morgan, *The Differential and Integral Calculus: Containing Differentiation, Integration, Development, Series, Differential Equations, Differences, Summation, Equations of Differences, Calculus of Variations, Definite Integrals* (London: Baldwin & Cradock, 1842). For further testimony on her mathematical work with De Morgan, see Huskey and Huskey, "Lady Lovelace and Charles Babbage," quotes p. 309: AAL to her mother "I go on most delightfully with Mr De Morgan. What can I ever do to repay him?"; AAL to Babbage "I am now studying the Finite Differences . . . And in this I have more particular interest, because I know it bears directly on some of your business"; and AAL to her mother "The Mathematics & Mr. De Morgan going on very well indeed. You would be much pleased to see the heap of papers of my writing."

[29] Toole, *Ada*, quote p. 173 (numbers of Bernoulli).

Frend, doubted the existence of negative numbers, De Morgan was a modern mathematician of the first rank. In his book *Trigonometry and Double Algebra* (1849) he presented a geometrical interpretation of complex numbers, those with real and imaginary parts.

In January 1844 De Morgan wrote a lengthy and detailed confidential letter to Ada's mother, Lady Byron, making an acute assessment of Ada's unusual facility with mathematics. "I never expressed to Lady Lovelace my opinion of her as a student in these matters," De Morgan began. "The power of thinking on these matters which Lady L[ovelace] has always shown from the beginning of my correspondence with her, has been something so utterly out of the common way for any beginner, man or woman, that this power must be duly considered by her friends . . . whether they should urge or check her obvious determination . . . to get beyond the present bounds of knowledge." De Morgan rated Ada favorably with Maria Agnesi, the Italian author of a pioneering calculus textbook (1748), and far more highly than Mary Somerville.[30]

During the same years as his correspondence with Ada, De Morgan also facilitated the mathematical work of George Boole, later author of the landmark *Investigation of the Laws of Thought* (1854) and today widely hailed as the father of Boolean algebra. In his *Treatise on the Calculus of Finite Differences*, Boole suggests a useful historical insight relating the character of mid-nineteenth century mathematics to the capabilities of Babbage's analytical engine.[31] Originally the Bernoulli numbers were discovered in 1712-13 (by the Swiss mathematician Jacob Bernoulli and by the Japanese mathematician Seki Kōwa), to aid in such calculations as the summation of powers ($1^n + 2^n + 3^n + 4^n$ . . .). With their assistance Bernoulli computed the sum of the first 1,000 integers, each raised to the *tenth* power, a immense number 91,409,924,241,424,243,424,241,924,242,500, in (as he claimed) "less than half of a quarter of an hour."[32]

Subsequent work by the Swiss and Scottish mathematicians Euler and Maclaurin connected the mathematics of integral calculus to the summation of polynomial expressions (which are essentially sums of powers each multiplied by some coefficient). Transcendental mathematical functions as well as integral calculus could be expressed by polynomial expressions; the Bernoulli numbers appeared as coefficients in some of these polynomial series.

---

[30] Velma R. Huskey and Harry D. Huskey, "Lady Lovelace and Charles Babbage," *Annals of the History of Computing* 2 no. 4 (1980): 299-329, De Morgan quoted p. 326; and Massimo Mazzotti, "Maria Gaetana Agnesi: Mathematics and the Making of the Catholic Enlightenment," *Isis* 92 no. 4 (2001): 657-683.

[31] Originally published in 1860, George Boole's *Treatise on the Calculus of Finite Differences* (New York: Dover, 1960) deals with Bernoulli numbers in chapter 6.

[32] Jacques [Jacob] Bernoulli, "On the 'Bernoulli numbers'," in David Eugene Smith, *A Source Book in Mathematics* (New York: McGraw-Hill, 1929), 85-90, quote p. 90. In modern notation, the Bernoulli numbers are defined as the coefficients ($B_n$) for the series expression for an exponential generating function, as follows (see <mathworld.wolfram.com/BernoulliNumber.html>)

Thus, by summing up the correct polynomial (using repeated multiplications, squaring, cubing, et seq.) Babbage's analytical engine could calculate transcendental mathematical functions (such as sine and cosine) as well as evaluate integral-calculus expressions, so long as they could be expanded into polynomial series. The Bernoulli numbers could be used to simplify the notation and computation of certain polynomial series, and thus were a powerful aid to computing. Bernoulli achieved his remarkable computation by transforming the extensive summation of powers ($1^{10} + 2^{10} + 3^{10} \ldots 1000^{10}$) into a straightforward seven-term polynomial equation using the Bernoulli numbers (up to $B_{10}$ in the series).[33]

Since Ada's calculus studies with De Morgan likely drew on the calculus textbook he was writing during these years, it is relevant to review its treatment of the Bernoulli numbers. De Morgan used the Bernoulli numbers in treating polynomial series expansions for $e^x$, tangent, and cotangent; in the calculus of operations; and in convergent series for definite integrals.[34] One exercise (page 307 §163) presents a *general expression* for directly computing a specific Bernoulli number, in terms of its predecessors.

$$B_{n+1} = -\frac{n+1}{2^{n+1}-1} \cdot \frac{1}{2+\Delta} \cdot 0^n = -\frac{n+1}{2^{n+1}-1} \left( \frac{0^n}{2} - \frac{\Delta 0^n}{4} + \frac{\Delta^2 0^n}{8} - \ldots \pm \frac{\Delta^n 0^n}{2^{n+1}} \right)$$

So for n = 7, the expression for $B_{n+1}$ or $B_8$ is as follows (where the fractional terms are computations on earlier instances in the series):

For the value of $B_n$ (n=7) we have

$$-\frac{8}{255} \left\{ -\frac{1}{4} + \frac{126}{8} - \frac{1806}{16} + \frac{8400}{32} - \frac{16800}{64} + \frac{15120}{128} - \frac{5040}{256} \right\} = -\frac{1}{30}.$$

So, with Ada Lovelace's interest in Babbage's machines, her mathematical studies with Babbage, Somerville, and De Morgan, and her relentless curiosity, she was surprisingly well exposed to the advanced mathematics of the period and had ample background and motivation to delve into the computations that appeared in the Notes to the Sketch.[35]

---

[33] Sum of the first 1000 integers raised to the tenth power = $1/11 x^{11} + B_1 x^{10} + 5B_2 x^9 + 30B_4 x^7 + 42B_6 x^5 + 15B_8 x^3 + B_{10} x$, where x = 1000 and $B_n$ are the Bernoulli numbers, viz., $B_1 = 1/2$, $B_2 = 1/6$, $B_4 = -1/30$, $B_6 = 1/42$, $B_8 = -1/30$, $B_{10} = 5/66$.

[34] See De Morgan, *The Differential and Integral Calculus*, pp. 247, 248, 308, 553, 581.

**Steps to the Sketch**

Babbage conceived a general computing machine around 1834, setting aside his still-uncompleted work on the Difference Engine to take up the challenges of what became known as the Analytical Engine. Whereas it was necessary to mechanically set up the Difference Engine to do each calculation, such as the prime-number generating polynomial $x^2 + x + 41$, Babbage's inspiration for the Analytical Engine was a computing machine able to re-configure itself—if not precisely "programmable" in the modern sense of the term.[36] Babbage's mature design while remaining a mechanical-age technology had several features in common with modern computers: separating the computation of numbers from their storage (he used the terms 'mill' and 'store' in an analogy with the industrial factory); adopting a mechanism to do conditional testing on intermediate results and hence permitting the branching of calculations; and using punched cards loosely inspired by the Jacquard loom.[37]

Babbage took one of his sets of evolving plans for the Analytical Engine to present at a conference in Turin in 1840. In the audience was a future prime minister of Italy. At the time Luigi Menabrea was a professor of mechanics and construction at the university of Turin; he subsequently served as a military engineer, naval minister, and eventually prime minister of Italy (1867-69).[38]  In October 1842, Menabrea published a short description of Babbage's Analytical Engine in a Swiss journal (written in French).  By this time, Babbage was well on the way to ruining whatever chance might have remained for support from the British government, especially after a disastrous meeting in November of that year with prime minister Robert Peel.  "What shall we do to get rid of Mr. Babbage and his calculating Machine?  Surely if completed it would be worthless as far as science is concerned," he wrote.  Peel, signaling his displeasure, soon dispatched Babbage's prime-number-calculating Difference Engine to the King's College Museum.[39]

---

[35] See also Imogen Forbes-Macphail's chapter in this volume [*Ada's Legacy*] exploring the "poetical" nature of mathematics.

[36] Babbage spent years working out a notation for expressing how its computations might be expressed. See Allan G. Bromley, "Charles Babbage's Analytical Engine, 1838," *IEEE Annals of the History of Computing* 20 no. 4 (1998): 29-45; and Allan G.Bromley, "Babbage's Analytical Engine Plans 28 and 28a: The Programmer's Interface," *IEEE Annals of the History of Computing* 22 no. 4 (2000): 5-19.

[37] The easy one-to-one correspondence between Jacquard looms and Babbage's cards, posited by such authors as James Essinger (author of popular works on Jacquard and Lovelace), is critically scrutinized by Martin Davis and Virginia Davis, "Mistaken Ancestry: The Jacquard and the Computer," *Textile* 3 no. 1 (2005): 76-87. After years of close study, Allan Bromley felt that the Analytical Engine fell somewhat short of a modern computer, as he wrote (privately) to Wilkes: "Perhaps my disappointment comes from being forced to accept that Babbage did NOT devise a COMPUTER but only a very sophisticated CALCULATOR after the style of the Harvard Mark I or the NCR accounting machines. He did not cross the watershed that marked off the [stored-program computers] EDVAC/EDSAC, although many of his technical innovations in implementation were astounding. It is a little difficult to admit this conclusion after expending so many years studying Babbage's designs. Perhaps it is why I 'took a break from Babbage' in the late 1980s and never really came back." See "Appreciating Charles Babbage: Emails between Allan Bromley and Maurice Wilkes," *IEEE Annals of the History of Computing* 26 no. 4 (2004): 62-70.

It was not initially Babbage who encouraged Ada Lovelace to examine the Menabrea manuscript and translate it into English.  Rather the prompting came from the notable scientist and sometime telegraph inventor Charles Wheatstone, who knew both Babbage and Lovelace and recruited contributions for *Taylor's Scientific Memoirs*.  Wheatstone had gained fame in 1837 with the patenting of a multiple-wire electric telegraph system using the positioning of magnetic needles to encode individual letters, rather than the Morse code system of dots and dashes.  He, too, was intrigued with using electromechanical apparatus for computation.  "Yesterday saw Wheatstone's model for telegraph and his drawings for Multiplication Engine," wrote Babbage after a visit in 1839.[40]  In 1843, the publication year of the Menabrea translation, Wheatstone improved on and publicized the famous "Wheatstone bridge" used to precisely measure electrical resistances.

Over the winter of 1842-43 Lovelace worked on translating the Menabrea manuscript, around 8,000 words, and first showed her results to Babbage in the spring of 1843.  In his memoir, *Passages from the Life of a Philosopher*, Babbage recalled encouraging Ada to add descriptive notes to her translation.

> We discussed together the various illustrations that might be introduced: I suggested several, but the selection was entirely her own.  So also was the algebraic working out of the different problems, except, indeed, that relating to the numbers of Bernoulli, which I had offered to do to save Lady Lovelace the trouble.  This she sent back to me for an amendment, having detected a grave mistake which I had made in the process.
>
> The notes of the Countess of Lovelace extend to about three times the length of the original memoir.  Their author has entered fully into almost all the very difficult and abstract questions connected with the subject.[41]

---

[38] See "Luigi Federico Menabrea" at <www-history.mcs.st-andrews.ac.uk/Biographies/Menabrea.html>.
[39] Fuegi and Francis, "Lovelace & Babbage," quote p. 16 (Peel on Babbage).
[40] Fuegi and Francis, "Lovelace & Babbage," quote p. 18 (Babbage on Wheatstone); Stein, *Ada*, 88.
[41] Charles Babbage, *Passages from the Life of a Philosopher* (London: Longman, Green, 1864), quote p. 136.

By this time in his life, while he was certainly writing for posterity and obviously keen on memorializing his computing engines, Babbage had no particular reason to exaggerate Ada's achievements. She had died years earlier, at age 36, and left Babbage a modest legacy of £600. I think we can take him at his word when he admits a "grave mistake" in deriving the Bernoulli numbers (which Bromley labels anachronistically as a "bug"). This passage, along with the Babbage-Lovelace letters, clearly describes a collaboration where Babbage and Lovelace are working together on the Notes. It's also clear from the context that "their author" here refers to Lovelace (rather than Menabrea) and that it is certain praise that *she* "has entered fully into . . . difficult and abstract questions." At the very least, Babbage's statement challenges the assertion that the notes "give an excellent account of his [that is, Babbage's alone] understanding of the mechanization of computational processes and the mathematical powers of the machine."

      Letters from the exact weeks in the summer of 1843 when Babbage and Lovelace were working on the Notes provide additional detail on their collaboration. This documentary evidence—from the British Library, the Science Museum, and the Oxford Bodleian Library—has been analyzed by Velma Huskey and Harry Huskey as well as John Fuegi and Jo Francis and published in *Annals of the History of Computing*.[42] The picture is clearly of a collaboration where both Babbage and Lovelace are making important contributions. Lovelace wrote the Notes and most of her letters to Babbage at her Ockham Park estate an hour's south of London; Babbage received her letters and responded from his London house on Dorset Street; and from time to time they met in person at Lovelace's London house on aristocratic St. James Square. Again, I emphasize that to point out Lovelace's knowledge, achievements, and contributions is *not* to belittle Babbage's.

      It is not easy to support the conjecture that the notes are Babbage's alone or that he directed Lovelace to write them. Bromley is not the only critic who has aimed to reduce Lovelace to a low-level clerk-assistant to Babbage. In *Ada: A Life and a Legacy*, Dorothy Stein, for instance, seized on a misprint that Lovelace, Babbage, and Menabrea before them all failed to spot—supposedly highlighting "the significance of her curiously ignored translation of a printer's error"—and contrives an argument that Lovelace had only a tenuous grasp of mathematics and a "dubious" understanding of the Analytical Engine's mechanical and logical operations. Stein contends, on this line of conjecture, that Lovelace was "completely dependent

---

[42] Velma R. Huskey and Harry D. Huskey, "Lady Lovelace and Charles Babbage," *Annals of the History of Computing* 2 no. 4 (1980): 299-329; and John Fuegi and Jo Francis, "Lovelace & Babbage and the Creation of the 1843 'Notes'," *IEEE Annals of the History of Computing* 25 no. 4 (2003): 16-26.

on [Babbage] for information and claims about the Analytical Engine."[43]  Stein's use of evidence is rather thin and highly selective.  In their analysis, Fuegi and Francis point to contemporaneous letters from Babbage, Lovelace, Wheatstone, the editors of *Taylor's Scientific Memoirs* and other scientific colleagues, concluding "all contemporaries of Lovelace and Babbage, having first-hand knowledge of how the 'Notes' came into being, acknowledged Lovelace at the time as the primary author."[44]

The evidence from Babbage's letters points to a collaboration between them, while their informal manner of addressing each other indicates a degree of collegiality.  Ada writes to "My Dear Babbage" while he responds to "My Dear Lady Lovelace."  On 30 June 1843 he writes:

> I am delighted with Note D.  It is in your usual clear style and required only one trifling alteration which I will make.  This arises from our not having yet had time to examine the outline of the mechanical part . . . I enclose a copy of the integration.  I am still working at some most entangled notations of Division but see my way through them at the expense of heavy labour . . . . Your latest information was the most agreeable.[45]

Let us turn to Note G on the computation of the Bernoulli numbers.  Recall that they were a common topic in 19th century mathematics, and that Lovelace herself had previously inquired about them to De Morgan.  The topic first appears in the correspondence when Ada wrote to Babbage, in a letter dated simply "Monday," as follows:

> I am working *very* hard for you . . . . I *think* you will be pleased.  I have made what appears to me some very important extensions & improvements . . . .
>
> I want to put something about Bernoulli's Numbers, in one of my Notes, as an example of how an implicit function may be worked out by the engine, without having been worked our by human head & hands first [as the Difference Engine required].  Give me the necessary data & formulae.[46]

Even if Babbage provided Ada with the mathematical expressions for the Bernoulli numbers, and assisted with the derivation of a general formula, the transformation of the general formula into a step-by-step algorithm remains Ada's achievement, as the letters clearly indicate.  The mathematics is somewhat involved, beginning with the basic equation where the Bernoulli numbers (note the *odd*-numbered notation $B_1$, $B_3$, $B_5$) appear as coefficients for the exponential function:

$$\frac{x}{\epsilon^x - 1} = 1 - \frac{x}{2} + B_1 \frac{x^2}{2} + B_3 \frac{x^4}{2 \cdot 3 \cdot 4} + B_5 \frac{x^6}{2 \cdot 3 \cdot 4 \cdot 5 \cdot 6} + \cdots$$

and then involving expansion, division, derivation, rounds of intricate multiplication, and finally writing the equation in general form:

$$\left. \begin{array}{l} 0 = -\frac{1}{2} \cdot \frac{2n-1}{2n+1} + B_1 \left(\frac{2n}{2}\right) + B_3 \left(\frac{2n \cdot (2n-1) \cdot (2n-2)}{2 \cdot 3 \cdot 4}\right) + \\ + B_5 \left(\frac{2n \cdot (2n-1)\ldots(2n-4)}{2 \cdot 3 \cdot 4 \cdot 5 \cdot 6}\right) + \cdots + B_{2n-1} \end{array} \right\}$$

---

[43] Stein, *Ada*, quote pp. xi (printer's error) 90-91 (tenuousness, dubious, completely dependent).  The printer's error was made in Menabrea's original publication where the French word *cas* (as in the case where N goes to infinity in a math expression) was mistakenly printed as *cos.* (as in the abbreviation for cosine).
[44] Fuegi and Francis, p. 26 note 29.
[45] Babbage to AAL 30 June 1843, quoted in Huskey and Huskey, pp. 312-313.
[46] AAL to Babbage, "Monday," quoted in Huskey and Huskey, pp. 311.

This equation (as note G explains, introducing the *odd*-numbered notation) "enables us to find . . . any nth Number of Bernoulli $B_{2n-1}$, in terms of all the preceding ones, if we but know the values of $B_1$, $B_3$ . . . $B_{2n-3}$." If n = 1, then the numerators for each of the higher terms ($B_3$ et seq.) contain a *zero* (2n-2) and so they drop out, permitting the direct calculation of $B_1$. Similarly, for n = 2 one substitutes the already computed value for $B_1$ while the higher terms ($B_5$ et seq.) again drop out and permit the direct computation of $B_3$. With n = 3, the process is repeated to compute $B_5$. The Notes explicitly make the general argument: "And so on, to any extent."

  Ada was determined to get the advanced mathematics correct. The two often exchanged letters daily, and sometimes even more frequently, as when a servant came into London and then waited for Babbage's reply. "I am doggedly attacking . . . all the ways of deducing the Bernoulli Numbers," she wrote at one moment. On the mathematics Ada wrote to Babbage (in a letter dated simply Tuesday morning) "the few lines I enclosed you last night about the connexion of (8) [the Bernoulli number equation in general form] with the famous Integral, I by no means intend you to insert, unless you *fully* approve the doing so."[47]

  The table and diagram that contained and expressed the algorithm were of special concern. Ada wrote Babbage (Saturday 6 o'clock), in connection with the note on the Bernoulli numbers, "Think of my horror then at just discovering that the Table & Diagram, (over which I have been spending infinite patience and pains) are seriously *wrong*, in one or two points. I have done them however in a beautiful manner, much improved upon our first edition of a Table and

---

[47] Toole, *Ada*, quote 204 (doggedly attacking); AAL to Babbage, "Tuesday morng," quoted in Huskey and Huskey, pp. 314.

Diagram. But unluckily I have made some errors. I send you this final note [G] excepting the Table & Diagrams." As a postscript, Ada adds: "Let me know how you like my finishing up of [G]. Mind you *scrutinise* all the *n*'s very carefully. I mean those of Sheets 4 and 5." Then, after a full day's work on Sunday, "you will admire the *Table & Diagram* extremely. They have been made out with extreme care & and all the *indices* most minutely & scrupulously attended to. Lord L[ovelace] is at this moment *inking it all over* for me. I had to do it in pencil."[48]

Ada tackled the question of the looping structure with some care. In a letter dated simply Tuesday, she writes to Babbage: "I hope you will approve of what I send. I have taken much pains with it. I have explained that there would be, in this instance & in many others, a recurring group or cycle of *Variable* as well as of *Operation* cards . . . ." She specifically identifies the main or outer loop as the "repetitions of (13 . . . 23)," and explains "as the variations follow a regular rule, they would be easily provided for."[49]

In his correspondence with Ada, Babbage directly admires her work on the Notes. His praise of Note D, "in your usual clear style," was cited earlier. In response to the Saturday letter quoted in the above paragraph, and clearly written before he received the Sunday letter, Babbage wrote: "I like much the improved form of the Bernoulli Note but can judge of it better when I have the diagram and Notation. I am very reluctant to return the admirable and philosophic view of the Anal. Engine contained in Note A. Pray do not alter it and do let me have it returned on Monday."[50] (Babbage was assembling the entire set of Notes in preparation for handing them to Charles Wheatstone for publication.)

As the final versions were going to the printer, Lovelace and Babbage were still revising the Notes' treatments of the variable cards and the operation cards, which specified the elemental arithmetical operators (+ - x ÷). Ada made clear that handling the operations cards was at the center of the Analytical Engine's capability for looping, even while the mechanical apparatus for effecting conditional tests was at the time not entirely clear.[51] Note D presents a straightforward computation—similar to one that Menabrea had presented in the original essay, also in tabular form—with six variable cards (for constants), nine working variables (for intermediate results), and two variables for results. There were 11 sequential steps (no looping structures), and the entirety fits onto the page (figure 1).

---

[48] Huskey and Huskey, pp. 313-14. In the earlier [Saturday] letter, Lovelace mistakenly labeled the final note "H" rather than "G," which was corrected in the later Sunday letter.

[49] Huskey and Huskey, pp. 316 (recurring group).

[50] Huskey and Huskey, pp. 313 (Babbage on improved form of the Bernoulli note).

[51] For a mechanical visualization of the Analytical Engine's "Logic and Loops," see Sydney Padua's *The Thrilling Adventures of Lovelace and Babbage*, pp. 306-308.

**Figure 1: Table-algorithm derived from Menabrea original (Note D)**

| Number of Operations | Nature of Operations | Variables for Data | | | | | | Working Variables | | | | | | | | | Variables for Results | |
|---|---|---|---|---|---|---|---|---|---|---|---|---|---|---|---|---|---|---|
| | | $^1V_0$ | $^1V_1$ | $^1V_2$ | $^1V_3$ | $^1V_4$ | $^1V_5$ | $^0V_6$ | $^0V_7$ | $^0V_8$ | $^0V_9$ | $^0V_{10}$ | $^0V_{11}$ | $^0V_{12}$ | $^0V_{13}$ | $^0V_{14}$ | $^0V_{15}$ | $^0V_{16}$ |
| | | + | + | + | + | + | + | + | + | + | + | + | + | + | + | + | + | + |
| | | 0 | 0 | 0 | 0 | 0 | 0 | 0 | 0 | 0 | 0 | 0 | 0 | 0 | 0 | 0 | 0 | 0 |
| | | 0 | 0 | 0 | 0 | 0 | 0 | 0 | 0 | 0 | 0 | 0 | 0 | 0 | 0 | 0 | 0 | 0 |
| | | 0 | 0 | 0 | 0 | 0 | 0 | 0 | 0 | 0 | 0 | 0 | 0 | 0 | 0 | 0 | 0 | 0 |
| | | 0 | 0 | 0 | 0 | 0 | 0 | 0 | 0 | 0 | 0 | 0 | 0 | 0 | 0 | 0 | 0 | 0 |
| | | | $m$ | $n$ | $d$ | $m'$ | $n'$ | $d'$ | | | | | | | | | | $\frac{dn'-d'n}{mn'-m'n}=x$ | $\frac{d'm-dm'}{mn'-m'n}=y$ |
| 1 | × | | $m$ | …. | …. | …. | $n'$ | …. | $mn'$ | | | | | | | | | | |
| 2 | × | | …. | $n$ | …. | $m'$ | …. | …. | …. | $m'n$ | | | | | | | | | |
| 3 | × | | …. | …. | $d$ | …. | …. | …. | …. | …. | $dn'$ | | | | | | | | |
| 4 | × | | …. | 0 | …. | …. | …. | $d'$ | …. | …. | …. | $d'n$ | | | | | | | |
| 5 | × | 0 | …. | …. | …. | …. | …. | 0 | …. | …. | …. | …. | $d'm$ | | | | | | |
| 6 | × | …. | …. | 0 | 0 | …. | …. | …. | …. | …. | …. | …. | $dm'$ | | | | | | |
| 7 | − | …. | …. | …. | …. | …. | …. | 0 | 0 | …. | …. | …. | …. | $(mn'-m'n)$ | | | | | |
| 8 | − | …. | …. | …. | …. | …. | …. | …. | …. | 0 | 0 | …. | …. | …. | $(dn'-d'n)$ | | | | |
| 9 | − | …. | …. | …. | …. | …. | …. | …. | …. | …. | …. | 0 | 0 | …. | …. | $(d'm-dm')$ | | | |
| 10 | ÷ | …. | …. | …. | …. | …. | …. | …. | …. | …. | …. | …. | …. | $(mn'-m'n)$ | 0 | …. | $\frac{dn'-d'n}{mn'-m'n}=x$ | |
| 11 | ÷ | …. | …. | …. | …. | …. | …. | …. | …. | …. | …. | …. | …. | 0 | …. | 0 | …. | $\frac{d'm-dm'}{mn'-m'n}=y$ |

For Note G on the Bernoulli numbers, the table–algorithm has ten data variables, three working variables, and four result variables. The computation has just 25 operations, but there are in addition two nested loops: an outer loop consisting of steps 13-23, and two inner loops consisting of steps 13-16 and 17-20. This form of calculation could not possibly have been completed with Babbage's Difference Engine, since it entirely lacked the ability to do conditional tests and create looping or branching structures. Nothing like it appeared in Menabrea's original. Lady Lovelace chose well when she identified the Bernoulli numbers as a means to show "how an implicit function may be worked out by the engine."

## Figure 2: Original table-algorithm for Bernoulli numbers (Note G)

Oversize table representing Ada Lovelace's algorithm for computing Bernoulli numbers. The nested looping structure—"here follows a repetition of operations thirteen to twenty-three"—is clearly visible at lower left (outer loop steps 13-23 and inner loops steps 13-16 and 17-20). Printed oversize and interleaved in *Taylor's Scientific Memoirs* (**high-res image** here **and** here).

This essay assesses Ada Lovelace's contribution to computing. It points out her mathematical studies with Babbage and De Morgan, her translation of Menabrea's Sketch, and her joint authorship with Babbage of the explanatory Notes. One element, the table-algorithm for computing the series of Bernoulli numbers, is by available evidence her work (written in pencil and inked in by her husband). Her correspondence with Babbage evinces a direct collegiality; the two figures were jointly grappling with how to communicate the details of a computing machine that did not physically exist. For this reason, the claim that Ada Lovelace was the world's first computer programmer might need to be carefully qualified. At the least, we can grant her primary authorship of the first algorithm intended for a computing machine.[52]

---

[52] See ACM's Collected Algorithms at <netlib.org/toms>.

It is surprising how widely and warmly she was recognized among early figures in computing. Alan Turing for instance felt obliged to deal with "Lady Lovelace's objection" to artificial intelligence since (in Turing's quotation) she had maintained "The Analytical Engine has no pretensions to originate anything. It can do whatever we know how to order it to perform."[53] A pioneering conference volume, *Faster Than Thought: A Symposium of Digital Computing Machines* (1953) was the source, according to the *History of Programming Languages* chapter on the Ada computer language, "from which we all learned about Lady Lovelace."[54] "Lady Lovelace had undoubtedly a profound understanding of the principles of the machine, and she added greatly to the value of her translation by some comprehensive notes about the machine . . . including what we should now call a *programme* for computing the Bernoulli numbers by a very sophisticated method," wrote a computer pioneer in 1953. The statement accords well with my assessment, after several decades of possibly fanciful writing about Ada Lovelace as well as unjustifiably critical blasts against her. Perhaps in the coming generation, we can come back to the measured appreciation of her achievements that originated, like the computer age itself, in the early 1950s. After all, as Ada Lovelace put it, "we may consider the [analytical] engine as the material and mechanical representative of analysis."[55]

---

[53] Andrew Hodges, "Turing: A Natural Philosopher," (1997) at https://www.turing.org.uk/publications/philobook.html ; and Darren Abramson, "Turing's Responses to Two Objections," *Minds and Machines* 18, no. 2 (2008): 147-167.

[54] Thomas J. Bergin, Jr. and Richard G. Gibson, Jr., eds., *History of Programming Languages---II* (New York ACM, 1996), "Ada Session," quote p. 207; B. V. Bowden, ed., *Faster Than Thought: A Symposium of Digital Computing Machines* (London: Isaac Pitman, 1953).

[55] B. V. Bowden, ed., *Faster Than Thought: A Symposium of Digital Computing Machines* (London: Isaac Pitman, 1953), 18-22, 70-75, 341-408, quote p. 18 (profound understanding); Lovelace, "Translator's Notes," quote p. 696 (representative of analysis). Compare Larry Owens, "Vannevar Bush and the Differential Analyzer: The Text and Context of an Early Computer," *Technology and Culture* 27 no. 1 (1986): 63-95.